\documentclass{conm-p-l}
\usepackage{amssymb,latexsym}
\usepackage{xypic}

\newcommand{\Mdef}[2]{\newcommand{#1}{\relax \ifmmode #2 \else $#2$\fi}}

\newcommand{\sm }{\wedge}

\newcommand{\tensor}{\otimes}

\newcommand{\ths}{\ltimes}        

\newcommand{\map}{\mathrm{map}}
\newcommand{\Hom}{\mathrm{Hom}}

\newcommand{\Ext}{\mathrm{Ext}}
\Mdef{\bhom}{\mathbf{\hat{H}om}}

\Mdef{\Mod}{\mathrm{mod}}

\newtheorem{thm}{Theorem}[section]
\newtheorem{lemma}[thm]{Lemma}
\newtheorem{prop}[thm]{Proposition}
\newtheorem{cor}[thm]{Corollary}

\theoremstyle{definition}

\newtheorem{defn}[thm]{Definition}

\newtheorem{construction}[thm]{Construction}

\newtheorem{example}[thm]{Example}

\newtheorem{remark}[thm]{Remark}

\newcommand{\qqed}{\qed \\[1ex]}

\Mdef{\PH} {\Phi^H}
\Mdef{\PK} {\Phi^K}
\Mdef{\PL} {\Phi^L}
\Mdef{\PT} {\Phi^{\T}}

\Mdef{\ef}{E{\cF}_+}
\Mdef{\etf}{\tilde{E}{\cF}}
\Mdef{\eg}{E{G}_+}
\Mdef{\etg}{\tilde{E}{G}}
\Mdef{\tf}{\T / \! \!  / {\cF}_+}

\Mdef{\infl}{\mathrm{inf}}
\Mdef{\defl}{\mathrm{def}}
\Mdef{\res}{\mathrm{res}}
\Mdef{\ind}{\mathrm{inf}}

\Mdef{\univ}{\mathcal{U}}

\Mdef{\Fp}{\mathbb{F}_p}
\Mdef{\Zpinfty}{\Z /p^{\infty}}
\Mdef{\Zpadic}{\Z_p^{\wedge}}

\newcommand{\bi}{\begin{itemize}}
\newcommand{\be}{\begin{enumerate}}
\newcommand{\bc}{\begin{center}}
\newcommand{\bd}{\begin{description}}
\newcommand{\ei}{\end{itemize}}
\newcommand{\ee}{\end{enumerate}}
\newcommand{\ec}{\end{center}}
\newcommand{\ed}{\end{description}}

\newcommand{\adjunction}[4]{
\diagram
#1:#2 \rrto<0.7ex> &&
#3  \llto<0.7ex> :#4 
\enddiagram}
\newcommand{\lra}{\longrightarrow}

\Mdef{\we}{\mathbf{we}}
\Mdef{\fib}{\mathbf{fib}}
\Mdef{\cof}{\mathbf{cof}}
\Mdef{\BI}{\mathcal{BI}}

\newcommand{\ilim}{\mathop{ \mathop{\mathrm{lim}} \limits_\leftarrow} \nolimits}
\newcommand{\colim}{\mathop{  \mathop{\mathrm {lim}} \limits_\rightarrow} \nolimits}
\newcommand{\holim}{\mathop{ \mathop{\mathrm {holim}} \limits_\leftarrow} \nolimits}

\Mdef{\A}{\mathbb{A}}
\Mdef{\B}{\mathbb{B}}
\Mdef{\C}{\mathbb{C}}
\Mdef{\D}{\mathbb{D}}
\Mdef{\E}{\mathbb{E}}
\Mdef{\T}{\mathbb{T}}
\Mdef{\F}{\mathbb{F}}
\Mdef{\G}{\mathbb{G}}
\Mdef{\I}{\mathbb{I}}
\Mdef{\N}{\mathbb{N}}
\Mdef{\Q}{\mathbb{Q}}
\Mdef{\R}{\mathbb{R}}
\Mdef{\bbS}{\mathbb{S}}
\Mdef{\Z}{\mathbb{Z}}

\Mdef{\bA}{\mathbb{A}}
\Mdef{\bB}{\mathbb{B}}
\Mdef{\bC}{\mathbb{C}}
\Mdef{\bD}{\mathbb{D}}
\Mdef{\bE}{\mathbb{E}}
\Mdef{\bF}{\mathbb{F}}
\Mdef{\bG}{\mathbb{G}}
\Mdef{\bH}{\mathbb{H}}
\Mdef{\bI}{\mathbb{I}}
\Mdef{\bJ}{\mathbb{J}}
\Mdef{\bK}{\mathbb{K}}
\Mdef{\bL}{\mathbb{L}}
\Mdef{\bM}{\mathbb{M}}
\Mdef{\bN}{\mathbb{N}}
\Mdef{\bO}{\mathbb{O}}
\Mdef{\bP}{\mathbb{P}}
\Mdef{\bQ}{\mathbb{Q}}
\Mdef{\bR}{\mathbb{R}}
\Mdef{\bS}{\mathbb{S}}
\Mdef{\bT}{\mathbb{T}}
\Mdef{\bU}{\mathbb{U}}
\Mdef{\bV}{\mathbb{V}}
\Mdef{\bW}{\mathbb{W}}
\Mdef{\bX}{\mathbb{X}}
\Mdef{\bY}{\mathbb{Y}}
\Mdef{\bZ}{\mathbb{Z}}

\Mdef{\cA}{\mathcal{A}}
\Mdef{\cB}{\mathcal{B}}
\Mdef{\cC}{\mathcal{C}}
\Mdef{\mcD}{\mathcal{D}} 
\Mdef{\cE}{\mathcal{E}}
\Mdef{\cF}{\mathcal{F}}
\Mdef{\cG}{\mathcal{G}}
\Mdef{\mcH}{\mathcal{H}} 
\Mdef{\cI}{\mathcal{I}}
\Mdef{\cJ}{\mathcal{J}}
\Mdef{\cK}{\mathcal{K}}
\Mdef{\mcL}{\mathcal{L}}

\Mdef{\cM}{\mathcal{M}}
\Mdef{\cN}{\mathcal{N}}
\Mdef{\cO}{\mathcal{O}}
\Mdef{\cP}{\mathcal{P}}
\Mdef{\cQ}{\mathcal{Q}}
\Mdef{\mcR}{\mathcal{R}}
\Mdef{\cS}{\mathcal{S}}
\Mdef{\cT}{\mathcal{T}}
\Mdef{\cU}{\mathcal{U}}
\Mdef{\cV}{\mathcal{V}}
\Mdef{\cW}{\mathcal{W}}
\Mdef{\cX}{\mathcal{X}}
\Mdef{\cY}{\mathcal{Y}}
\Mdef{\cZ}{\mathcal{Z}}

\Mdef{\tA}{\tilde{A}}
\Mdef{\tB}{\tilde{B}}
\Mdef{\tC}{\tilde{C}}
\Mdef{\tE}{\tilde{E}}
\Mdef{\tH}{\tilde{H}}
\Mdef{\tK}{\tilde{K}}
\Mdef{\tL}{\tilde{L}}
\Mdef{\tM}{\tilde{M}}
\Mdef{\tN}{\tilde{N}}
\Mdef{\tP}{\tilde{P}}

\Mdef{\Ab}{\overline{A}}
\Mdef{\Bb}{\overline{B}}
\Mdef{\Cb}{\overline{C}}
\Mdef{\Db}{\overline{D}}
\Mdef{\Hb}{\overline{H}}
\Mdef{\Ib}{\overline{I}}
\Mdef{\Kb}{\overline{K}}
\Mdef{\Lb}{\overline{L}}
\Mdef{\Mb}{\overline{M}}
\Mdef{\Nb}{\overline{N}}
\Mdef{\Qb}{\overline{Q}}
\Mdef{\Tb}{\overline{T}}

\Mdef{\qb}{\overline{q}}
\Mdef{\rb}{\overline{r}}
\Mdef{\tb}{\overline{t}}
\Mdef{\ub}{\overline{u}}
\Mdef{\vb}{\overline{v}}
\Mdef{\hc}{\hat{c}}
\Mdef{\he}{\hat{e}}
\Mdef{\hf}{\hat{f}}
\Mdef{\hA}{\hat{A}}
\Mdef{\hH}{\hat{H}}
\Mdef{\hJ}{\hat{J}}
\Mdef{\hM}{\hat{M}}
\Mdef{\hP}{\hat{P}}
\Mdef{\hQ}{\hat{Q}}

\Mdef{\bolda}{\mathbf{a}}
\Mdef{\boldb}{\mathbf{b}}
\Mdef{\boldD}{\mathbf{D}}

\Mdef{\fm}{\frak{m}}

\Mdef{\eps}{\epsilon}

\newcommand{\extsmash}{\overline{\sm}}

\begin{document}

\newcommand{\po}{{\partial}}
\newcommand{\ov}{\overline}
\newcommand{\ga}{{\gamma}}
\newcommand{\uc}{{\mathcal U}}
\newcommand{\vc}{{\mathcal V}}
\newcommand{\wc}{{\mathcal W}}
\newcommand{\lc}{{\mathcal L}}
\newcommand{\lb}{{\mathbb{L}}}
\newcommand{\oc}{{\mathcal O}}
\newcommand{\re}{{\mathbb{R}}}
\newcommand{\un}{\text{unst}}
\newcommand{\basedspaces}{\mathbf{Spaces_*}}
\newcommand{\smbar}{\overline{\sm}}
\newcommand{\coeq}{\mathrm{coeq}}
\newcommand{\Rmod}{\mbox{$R$-}\mathrm{mod}}
\newcommand{\bRmod}{\mbox{$\boldR$-}\mathrm{mod}}
\newcommand{\modcE}{\mathrm{mod}\mbox{-$\cE$}}
\newcommand{\HombR}{\mathrm{Hom}_{\boldR}}
\newcommand{\Hombk}{\mathrm{Hom}_{\boldk}}
\newcommand{\tensorcE}{\otimes_{\cE}}
\newcommand{\HomcE}{\mathrm{Hom}_R}
\newcommand{\HomR}{\mathrm{Hom}_R}
\newcommand{\Homk}{\mathrm{Hom}_k}
\newcommand{\boldR}{\mathbf{R}}
\newcommand{\boldk}{\mathbf{k}}

\title{Spectra for commutative algebraists.}
\author{J.P.C.Greenlees}
\address{Department of Pure Mathematics, Hicks Building,
Sheffield S3 7RH. UK.}
\email{j.greenlees@sheffield.ac.uk}
\thanks{This article grew out of a series of talks given as part of the 
MSRI emphasis year on commutative algebra. JPCG is grateful to L.Avramov 
for the invitation.}
\date{}
\begin{abstract}
The article  is designed to explain to commutative algebraists what spectra 
are, why they were originally defined, and how they can be useful for commutative
 algebra. 
\end{abstract}
\maketitle

\tableofcontents

\section{Introduction.}
This article grew out of a short series of talks given as part of the 
MSRI emphasis year  on  commutative algebra. The purpose is to explain 
to commutative algebraists what spectra (in the sense of homotopy 
theory) are, why they were originally defined, and how they can be useful 
for commutative algebra. An account focusing on applications in commutative 
algebra rather than foundations can be found in \cite{Guanajuato}, and
an introduction to the methods of proof can be found in another article
in the present volume \cite{firststeps}. 

Historically, it was only after several refinements that spectra sufficiently rigid for 
the  algebraic applications were defined. We will follow a similar path, 
 so it will take some time before algebraic examples can be explained. 
Accordingly, we begin with an overview to explain where we are going. 
We only intend to give an outline and overview, not a course in homotopy 
theory: detail will be at a minimum, but we give references at the 
appropriate points for those who wish to pursue the subject further. 
As general background references we suggest \cite{May} for general 
homotopy theory leading towards spectra, \cite{GoerssJardine} for 
simplicial homotopy theory and \cite{Hovey} for Quillen model categories.
A very approachable introduction to spectra is given in \cite{Omega}, but most of 
the applications to commutative algebra have only become possible because of 
developments since it was written. The main foundational sources for 
spectra are collected
at the start of the bibliography, with letters rather than numbers for
their citations.

\section{Motivation via the derived category.}

Traditional commutative algebra considers commutative rings $R$ and
modules over them, but some constructions make it natural to extend
further  to considering chain complexes of $R$-modules; the need to 
consider robust, homotopy invariant properties leads to the derived 
category $D(R)$. Once we admit chain complexes, it is natural to 
consider the corresponding multiplicative objects, differential graded algebras.
Although it may appear inevitable, the real justification for this 
process of generalization is the array  of naturally occurring examples. 

The use of spectra is a logical extension of this process: 
they allow us to define  flexible generalizations of the derived category.
Ring spectra extend the notion of rings, module spectra extend the notion
of chain complexes, and the homotopy category of module spectra extends the
 derived category. 
Many ring theoretic constructions extend to ring spectra, 
and thus extend the power of commutative algebra to a vast new supply of naturally
occurring examples. Even for traditional rings, the new perspective is
often enlightening, and thinking in terms of spectra makes
 a number of new tools available. Once again the only compelling justification 
for this inexorable process of generalization  is the array of naturally occurring 
examples, some of which we will be described later in this article. 

We now rehearse some of the familiar arguments for the derived category of a ring in more 
detail, so that it can serve  as a model for the case of ring spectra. 
\vskip .1 in

\subsection{Why consider the derived  category?}
The category of $R$-modules has a lot of structure, but it 
is rather rigid, and not well designed for dealing with 
homological invariants and derived functors. The derived 
category $D(R)$ is designed for working with homological invariants
and other properties which are homotopy invariant: it inherits
structure from the module category, but in an adapted form.

\begin{description}
\item[Modules] Conventional $R$-modules give objects of the derived category. 
It therefore contains many familiar objects. On the other hand, it contains many 
other objects (chain complexes), but  all objects of the
derived category are constructed from modules.  
\item[Homological invariants] Tor, Ext, local cohomology  and other 
homological invariants are represented in $D(R)$ and the derived category
$D(R)$ provides a flexible  environment for manipulating them.
Indeed, one may view the derived category as the universal domain 
for homological invariants. After the construction of the derived 
category, homological invariants reappear as pale shadows of the objects 
which represent them.

This is one reason for including so many new objects in the derived category.
Because the homological invariants are now embodied, they may be 
very conveniently compared and manipulated.

\end{description}
The derived category inherits a lot of useful structure from the 
category of modules. 
\begin{description}
\item[Triangulation] In the abelian category of $R$-modules, 
kernels, cokernels and exact sequences allow one to measure 
how close a map is to an isomorphism. Passing to the derived category, 
short exact sequences give triangles, and the use of triangles gives a 
way to internalize the deviation from isomorphism.
\item[Sums, products] We work with the unbounded derived category
and therefore have all sums and products. 
\item[Homotopy direct and inverse limits] In the module category 
it is useful to be able to construct direct and inverse
limits of diagrams of modules. However these are not homotopy 
invariant constructions: if one varies the diagram by a homotopy, 
the resulting limit need not be homotopy equivalent to the 
original one.

The counterparts in the derived category are 
homotopy direct and inverse limits. Perhaps the most familiar case
is that of a sequence
$$\cdots \lra X_{n-1}\stackrel{f_{n-1}}\lra X_n \stackrel{f_n}\lra X_{n+1}
\lra \cdots .$$
One may construct the direct limit as the cokernel of the
map
$$
1-f : \bigoplus_n X_n \lra \bigoplus_n X_n.  
$$
The {\em homotopy} direct limit is the next term in the triangle
(the mapping cone of $1-f$).
Because the direct limit over a sequence is exact,  the 
construction is homotopy invariant, and the direct limit itself
provides a model for the homotopy direct limit.
Similarly one may construct the inverse limit as the kernel of the
map
$$
1-f : \prod_n X_n \lra \prod_n X_n, 
$$
and in fact the cokernel is the first right derived functor of the
inverse limit.
The {\em homotopy} inverse limit is the previous term in the triangle
(the mapping fibre of $1-f$). Because the inverse  limit functor is 
not usually exact, one obtains a short exact sequence
$$
0 \lra \ilim^1_nH_{i+1}(X_n) \lra H_i(\holim_n X_n)
\lra \ilim_n H_i(X_n) \lra 0.
$$

One useful example is that this allows one to split all idempotents. 
Thus if $e$ is an idempotent self-map of $X$, the corresponding summand 
is both the homotopy direct limit and the homotopy inverse limit 
of the sequence 
$(\cdots\lra X \stackrel{e}\lra X \stackrel{e}\lra X \lra \cdots )$.
\end{description}
\vskip .3in

\subsection{How to construct the derived category.}

The steps in  the construction of  the derived category $D(R)$ of a ring
or differential graded (DG) ring $R$ may 
be described as follows. We adopt a somewhat elaborate
approach so that it provides a template for the corresponding process 
for spectra.\\[1ex]

{\bf Step 0:} {\em Start with graded sets with cartesian product. }
This provides the basic environment within which the rest of the construction takes
place. However, we need to move to an additive category.

{\bf Step 1:} {\em Form the category of graded abelian groups.} 
This provides a more convenient and
algebraic environment. Next we need additional multiplicative structure.

{\bf Step 2:} {\em Construct and exploit the tensor product.}
\begin{description}
\item[Step 2a] Construct the tensor product $\otimes_{\mathbb Z}$.

\item[Step 2b] Define differential graded (DG) abelian groups.

\item[Step 2c] Find the DG-abelian group $\mathbb Z$ and recognize DG-abelian 
groups as DG-$\mathbb Z$-modules.
\end{description}

{\bf Step 3:} {\em Form the categories of differential graded rings and modules.}
First we take a DG-$\mathbb Z$-module $R$ with the structure of a ring 
in the category of DG-$\mathbb Z$-modules,  and then define modules over it.
This constructs the algebraic objects behind the derived category. Finally, 
we pass to homotopy invariant structures.

{\bf Step 4:} {\em Invert homology isomorphisms.} In one sense this is a purely categorical 
process, but to avoid set theoretic difficulties and to make it accessible
to computation, we need to {\em construct} the category with homology isomorphisms 
inverted. One way to do this is to restrict to complexes of $R$-modules which are projective
in a suitable sense, and then pass to homotopy; a  flexible language for expressing this
 is that of model categories. \\[2ex]

We may summarize this process in the picture 
$$\begin{array}{clcl}
\, & \text{(0)\,\, Graded sets} &&\\
& \qquad \downarrow &&\\
& \text{(1)\,\, $\mathbb Z$-modules}&& \\
&\qquad \downarrow &&\\
& \text{(2)\,\, $DG$-$\mathbb Z$-modules} &\longrightarrow 
                             &\text{(4) Ho}(\mathbb Z\text{-mod})=D(\mathbb Z) \\
&\qquad \downarrow &&\\
& \text{(3)\,\, $R$-modules} &\longrightarrow
                             &\text{(4) Ho}(R\text{-mod})=D(R) 
\end{array}$$

One of the things to note about this algebraic situation is that there is no 
direct route from the derived category $D(\Z)$ of $\Z$-modules to the derived 
category $D(R)$ of $R$-modules. 
We need $R$ to be an actual DG ring (rather than a ring object in $D(\Z)$), 
and to consider actual $R$-modules (rather 
than module objects in $D(\Z)$). The technical difficulties of this step for 
spectra took several decades to overcome.

The rest of the article will sketch how to parallel this development for 
spectra, with the ring $R$ replaced by a ring {\em spectrum} and modules 
over $R$ replaced by module {\em spectra} over the ring spectrum. 
First we give a very brief motivation for considering 
spectra in the first place, and it will not be until
 Section \ref{sec:smash} that it  becomes possible to 
explain what we mean by ring spectra. For the present we speak very 
informally,  not starting to give definitions until Section \ref{sec:prespectra}.

\section{Why consider spectra?} 

We will answer the question from the point of view of an algebraic
topologist. 

To avoid changing later, all our spaces will come equipped
with specified basepoints. We write $[X,Y]_{\un}$ for the set of 
based homotopy classes of based maps from $X$ to $Y$ and we write
$H^*(X)$ for the reduced cohomology of $X$ with integer coefficients. 
The subscript $\un$ is short for `unstable'; this is to contrast with 
`stable' maps of spectra, described below.

\subsection{First Answer.} Spectra describe a relatively well behaved part of
homotopy theory \cite{SpanierWhitehead}. We will see later that 
spaces give rise to spectra and, for highly connected
spaces,  homotopy classes of maps of spaces and of the corresponding
spectra coincide. 

To be more precise, we need the suspension functor $\Sigma Y := Y \sm S^1 $ 
where the smash product of based spaces  is 
$X \sm Y:= X \times Y/(X\times\{y_0\} \cup \{ x_0\} \times Y)$. If $X$ is 
a CW-complex,  the suspension $\Sigma X$ is a CW-complex with  cells 
corresponding to those of $X$, but one dimension higher. Now we   
define the morphisms in the {\em Spanier-Whitehead category} by
$$
[X,Y]:=\colim_{ k}[\Sigma^kX,\Sigma^kY]_{\un}, 
$$ 
where the limit is over the suspension maps
$$\Sigma : [\Sigma^kX,\Sigma^kY]_{\un} \lra [\Sigma^{k+1}X,\Sigma^{k+1}Y]_{\un}. $$
An element of this direct limit is called a  `stable' map.
In fact $[X,Y]$ is an abelian group,  because the first suspension coordinate 
allows addition by concatenation, and the second suspension coordinate 
gives room to move the terms added past each other. 
It will transpire that  when $X$ is finite dimensional
the group $[X,Y]$ gives the maps
from  the spectrum associated to $X$  to the spectrum associated to $Y$.
Furthermore, it turns out that the above limit is achieved, and hence the 
maps of spectra give a very well behaved piece of homotopy theory. 
To explain this, write $bottom (Y)$ for the lowest dimension of a cell 
in $Y$ and $dim (Y)$ for the highest. The Freudenthal suspension theorem  
states that suspension gives an isomorphism
$$
\Sigma: [X,Y]_{\un} \overset{\cong}{\longrightarrow} 
[\Sigma X,\Sigma Y]_{\un} 
\text{ if } \dim X\le 2 \cdot \text{ bottom}(Y)-2 .
$$
Thus if $X$ is finite dimensional  all the maps in the direct limit system
are eventually isomorphic.

One reason for considering stable maps is that the suspension isomorphism
$$
 H^n(X) \cong  H^{n+1}(\Sigma X) \cong  H^{n+2}(\Sigma^2 X) \cong \dots
$$
for reduced cohomology shows that  it is stable maps that are relevant
to cohomology. More precisely, if a stable map $f:X \lra Y$ is represented
by a continuous function $g: \Sigma^k X \lra \Sigma^k Y$, then $f$ induces a map 
$f^*$ in cohomology so that the diagram
$$\begin{array}{ccc}
H^n(Y) & \stackrel{f^*}\lra &H^n(X)\\
\cong \downarrow &&\downarrow \cong \\
H^{n+k}(\Sigma^k Y) & \stackrel{g^*}\lra &H^{n+k}(\Sigma^kX)\\
\end{array}$$
commutes.

\subsection{Second answer.} Cohomological invariants are represented. Indeed 
(Brown representability \cite{Brown}) any contravariant homotopy functor 
$E^*(\cdot )$ on spaces
which satisfies the Eilenberg-Steenrod axioms 
(Homotopy, Excision/Suspension, Exactness) and the wedge axiom,  is 
represented by a spectrum $E$ in the sense that
$$
 E^*(X)=[X,E]^*.
$$
This equation introduces Adams's convenient abbreviation whereby the
name of the functor  $E^*(\cdot) $ on the left has been used to provide
the name for the representing spectrum $E$ on the right. The convention 
is also used in the reverse direction to name the cohomology 
theory represented by a spectrum which already has a name.

In effect,  this gives a way of {\em constructing} spectra, and hence a
way of arguing geometrically with cohomology theories. 
For example the Yoneda lemma shows that natural transformations of 
cohomology theories which commute with suspension ({\em stable cohomology 
operations}) are represented: 
$$\text{Stable cohomology operations}(E^*(\cdot ), F^*(\cdot))=[E,F]^*.$$
In particular the stable  operations between $E^*(\cdot)$ and itself
form the ring  $E^*E=[E,E]^*$.

\subsection{Third answer.} Naturally occurring invariants occur as homotopy 
groups of spectra. For example various sorts of bordism, and
algebraic $K$-theory. Similarly, many invariants in geometric topology are
defined as homotopy groups of  classifying spaces, and very often these spaces
are the infinite loop spaces associated to spectra.  This applies to 
Quillen's algebraic $K$-groups, originally defined as 
the homotopy groups of the space $BGL(R)^+$:
there is a spectrum $K(R)$ with $K_*(R)=\pi_*(K(R))$. Examples from geometric 
topology include the Whitehead space $Wh(X)$, Waldhausen's $K$-theory of spaces  
$A(X)$ \cite{Waldhausen} and 
the classifying space of the stable mapping class group $B\Gamma_\infty^+$
\cite{Tillmann}.
We will give further details of some of these constructions later. 

\subsection{Fourth answer.} This, finally, is  relevant to 
commutative algebraists. Many of the invariants described above are not just
groups, but also rings. In many cases this additional structure is reflected 
geometrically in the sense that the representing spectra have a product making 
them into rings (or  even commutative rings) 
in a suitable category of spectra. 
These spectra {\it with an appropriate tensor product\/} provide a 
context like the derived category.

 Several familiar algebraic constructions on rings can then 
also be applied to ring spectra to give new spectra. For example 
Hochschild homology and cohomology extends to topological Hochschild 
homology and cohomology,  Andr\'e-Quillen
cohomology of commutative rings extends to topological  Andr\'e-Quillen 
cohomology of commutative ring spectra, and algebraic 
$K$-theory of rings extends to $K$-theory of ring spectra. 
We will give further details of some of these constructions later. 

\section{How to construct spectra (Step 1).}
\label{sec:prespectra}

The counterpart to the use of graded sets in Step 0 of algebra 
is the use of based spaces. This section deals with the  Step 1 
transition to an additive category (corresponding to the formation of abelian groups
in the algebraic case). Based on the discussion of the Freudenthal suspension theorem,  the
definition of a spectrum is fairly natural. For the present, we take 
a fairly naive approach, perhaps best reflected in \cite{shagh}, although
the approach in the first few sections of \cite{LMS(M)} is more appropriate
for later developments. 

We begin with the first approximation to a spectrum.  

\begin{defn}
\label{defn:spectra}
A {\it spectrum\/} $E$ is a sequence of based spaces $E_k$ for  $k\ge 0$ 
together with structure maps
$$
\sigma\colon \Sigma E_k \to E_{k+1}.
$$
A map of spectra $f\colon E\to F$ is a sequence of maps so that the
squares
$$
\begin{array}{ccc}
\Sigma E_k  &\stackrel{\Sigma f_k}\longrightarrow &  \Sigma F_k \\
\downarrow && \downarrow \\
E_{k+1}    &\stackrel{f_{k+1}}\longrightarrow   &   F_{k+1}
\end{array}
$$
commute for all $k$.
\end{defn}

\begin{remark} May and others would call this a `prespectrum', 
reserving `spectrum'  for the best sort of prespectrum. To 
avoid conflicts, we will instead add adjectives to restrict 
the type of spectrum.
\end{remark}

\begin{example}
If $X$ is a based space we may define the suspension spectrum 
$\Sigma^\infty X$ to have $k$th term $\Sigma^k X$ with the structure maps
being the identity.
\vskip .1in

\noindent
{\bf Remark:}
It is possible to make a definition of homotopy immediately, but this does not 
work very well for arbitrary spectra. Nonetheless it will turn out that for 
finite CW-complexes $K$, maps out of a suspension  spectrum are given by 
$$
[\Sigma^\infty K,E]=\colim_{ k} [\Sigma^kK,E_k]_{\un}.
$$
In particular
$$
\pi_n(E):=[\Sigma^\infty S^n,E]=\colim_{ k} [S^{n+k},E_k]_{\un}.
$$
For example if $E=\Sigma^\infty L$ for a based space $L$, 
we obtain the {\em stable homotopy groups}
$$
\pi_n(\Sigma^\infty L)=\colim_{ k} [\Sigma^kS^n,\Sigma^kL]_{\un}, 
$$
which coincides with the group of maps $[S^n,L]$ in the Spanier-Whitehead 
category. By the Freudenthal suspension theorem, this is the common 
stable value of the groups  $[\Sigma^kS^n,\Sigma^kL]_{\un}$ for large $k$.
Thus spectra have captured stable homotopy groups.
\end{example}

\begin{construction}
We can suspend spectra by any integer $r$, defining $\Sigma^rE$ by
$$
(\Sigma^rE)_k= \begin{cases}
E_{k-r} \quad & k-r\ge0 \\
pt \quad & k-r <0.
\end{cases}
$$

Notice that if we ignore the first few terms,  $\Sigma^r$ is inverse to $\Sigma^{-r}$.
Homotopy groups involve a direct limit and therefore do not see these first few
terms. Accordingly,   once we invert homotopy isomorphisms,  the suspension functor becomes
an equivalence of categories. Because suspension is an equivalence,  we say that 
 we have a {\em stable} category.

\end{construction}

\vskip .1in

\begin{example} In particular we have sphere spectra. We write
$\bbS =\Sigma^{\infty}S^0$ for the 0-sphere because of its special role, and
then define
$$
S^r=\Sigma^r \bbS \quad \text{ for all integers } r .
$$
Note that $S^r$ now has meaning for a space and a spectrum for $r \geq 0$, 
but since we have an isomorphism $S^r\cong \Sigma^\infty S^r$ of spectra
for $r\ge0$ the ambiguity is not important.
We extend this ambiguity, by often suppressing $\Sigma^\infty$.
\end{example}

\begin{example}
{\it Eilenberg-MacLane spectra\/}.
An Eilenberg-MacLane  space of type $(R,k)$ for a group $R$ and $k \geq 0$
is  a CW-complex $K(R,k)$ with $\pi_k(K(R,k))=R$ and $\pi_n(K(R,k))=0$ for 
$n\neq k$; any two such spaces are homotopy equivalent.
It is well known that each cohomology group is represented by an Eilenberg-MacLane space. Indeed, for any CW-complex $X$,  we have 
$ H^k(X;R)=[X,K(R,k)]_{\un}$. In fact,   this sequence of spaces, as $k$ varies,  
assembles to make a spectrum.

To describe this, first note that the suspension functor $\Sigma$  is defined
 by smashing with the circle $S^1$, so it is left adjoint to the loop 
functor $\Omega$ defined by $\Omega X := \map (S^1,X)$ (based loops, 
with a suitable topology). In fact there is a homeomorphism
$$
\map(\Sigma W,X)=\map (W \sm S^1, X) \cong \map (W , \map (S^1,X))=\map (W , \Omega X)
$$
 This passes to homotopy, so  looping shifts homotopy
in the sense that $\pi_n(\Omega X)=\pi_{n+1}(X)$. We conclude that there
is a homotopy equivalence
$$
\tilde{\sigma}: K(R,k)\overset{\simeq}{\to} \Omega K(R,k+1) , 
$$
and hence we may obtain a spectrum 
$$
HR=\{K(R,k)\}_{k \geq 0} 
$$
where the bonding map 
$$
\sigma\colon\Sigma K(R,k)\to K(R,k+1) 
$$
is adjoint to $\tilde\sigma$. Thus we find
$$
[\Sigma^r \Sigma^\infty  X, HR]= 
\colim_k [\Sigma^r\Sigma^kX, K(R,k)]_{\un}=
\colim_k  H^k(\Sigma^r \Sigma^kX;R)= H^{-r}(X;R).
$$
In particular the Eilenberg-MacLane spectrum has homotopy in a single degree like the 
spaces from which it was built: 
$$\pi_k(HR)= \begin{cases}
R \quad k=0 \\
0 \quad k\ne 0. \end{cases}$$
\end{example}

\begin{example} The classification of smooth compact manifolds provided 
an important motivation for the construction of spectra.  Although this may 
seem too geometric for applications to commutative algebra, rather 
 mysteriously the spectra that arise this way are amongst those with the
most algebraic behaviour.

If we consider two $n$-manifolds
to be equivalent if they together form the boundary of an $(n+1)$-manifold
(they are `cobordant') we obtain the set $\Omega_n^O$ of cobordism classes
of $n$-manifolds. 
The superscript $O$ stands for `orthogonal', and refers to 
the fact that a bundle over a manifold admits a Riemannian metric and hence 
the normal bundle of an $n$-manifold embedded in Euclidean space
has a reduction to the orthogonal group. The set $\Omega_n^O$
 is a group under disjoint union, and taking all $n$
together we obtain a graded commutative ring with product induced by cartesian
product of manifolds. The group $\Omega^O_n$ may be calculated as the $n$th homotopy 
group of a spectrum $MO$.
The idea is that a manifold $M^n$ is determined up to cobordism by specifying
an embedding in $\R^{N+n}$ and considering its normal bundle $\nu$. 
Collapsing the complement of the normal bundle defines the so-called 
{\em Thom space} $M^{\nu}$ of $\nu$ and the {\em Pontrjagin-Thom collapse
map} $S^{N+n} \lra M^{\nu}$. On the other hand,  the
normal bundle is $N$-dimensional and thus classified by a map 
$\xi_\nu : M \lra BO(N)$, where $BO(N)$ is the classifying space for 
$O(N)$-bundles with universal bundle $\gamma_N$ over it. Taking the Thom 
spaces and composing with the  collapse map,  we have
$$
S^{N+n} \lra M^{\nu } \lra BO(N)^{\gamma_N}.
$$
By embedding $\R^{N+n}$ in $\R^{N+N'+n}$ these maps for different $N$ may
be compared, and as $N$ gets large, the resulting class in 
$$
\colim_N[S^{N+n}, BO(N)^{\gamma_N}]_{\un}
$$
 is independent of the 
embedding, and only depends on the cobordism class of $M$. Furthermore,  the 
manifold $M$ can be recovered up to cobordism by taking the transverse 
inverse image of the zero section. This motivates the definition of the 
cobordism spectrum $MO$.

We take $MO(n) : = BO(n)^{\gamma_n} $
and the bonding map is 
$$
\Sigma MO(n) = BO(n)^{\gamma_n\oplus 1}=BO(n)^{i^*\gamma_{n+1}} 
\to BO(n+1)^{\gamma_{n+1}}=MO(n+1). 
$$
The motivating discussion of the Pontrjagin-Thom construction thus proves
$$
\pi_nMO = \colim_{N} [S^{n+N},MO(N)]_{\un} \cong \Omega_n^O.
$$
It is by this means that Thom calculated the group $\Omega_n^O$ of cobordism 
classes of $n$-manifolds \cite{Thom}.

There are many variants of this depending on the additional structure on 
the manifold. Of particular importance are manifolds with a complex
structure on their stable normal bundle. The group of bordism classes of
these is $\Omega^U_*$ (the superscript now refers to the fact that the 
stable normal bundle has a reduction to a unitary group), and again this 
is given by the  homotopy groups of the Thom spectrum $MU$, and this allowed
Milnor to calculate the complex cobordism ring
$$
\Omega^U_*=\pi_*MU=\Z [x_1, x_2, \ldots ]
$$
where $x_i$ has degree $2i$ \cite{MilnorMU}. The spectrum $MU$ plays a central role 
in stable homotopy, both conceptually and computationally. It provides a close 
link with various bits of algebra, and in particular with commutative
algebra. The root of this connection is Quillen's theorem \cite{Quillen}  that the polynomial 
ring is  isomorphic to  Lazard's universal ring for one dimensional commutative 
formal group laws for geometric reasons.
\end{example}

\begin{example}
The theory of vector bundles gives rise to topological $K$-theory. Indeed, 
the unreduced complex $K$-theory of an unbased compact space $X$ is given by 
$$
K(X)=Gr(\C\text{-bundles over } X ), 
$$
where $Gr$ is the Grothendieck group completion. 
The reduced theory is defined by $K^0(X)=\ker (K(X) \lra K(pt))$, and 
represented by the space $BU \times \Z$ in the sense that
$$
K^0(X)=[X,BU \times \Z]_{\un}.
$$
The suspension isomorphism allows one to define $K^{-n}(X)$ for $n \geq 0$, 
but to give $K^n(X)$ we need Bott periodicity \cite{Bott, AB}. In terms of the cohomology 
theory,  Bott periodicity states $K^{i+2}(X)\cong K^i(X)$, and in terms of
representing spaces it states
$$
\Omega^2(BU\times \mathbb Z)\simeq BU\times \mathbb Z.
$$
Hence we may define the representing spectrum  $K$ by giving it $2n$th 
term $BU\times \mathbb Z$ and 2-fold bonding maps adjoint to the 
Bott periodicity equivalence
$BU\times\mathbb Z \overset{\simeq}{\to}\Omega^2 (BU\times \mathbb Z)$. We then 
find
$$
[\Sigma^\infty X,K]=\colim_{ k}[\Sigma^{2k}X, BU\times\mathbb Z]_{\un}=
[X,BU\times\mathbb Z]_{\un}=K^0(X)
$$
\end{example}

\begin{remark} (a) Spectra  with the property 
 $\Omega E_{k+1}\simeq E_k$  for all $k$ are  called 
{\em $\Omega$-spectra} (sometimes pronounced `loop spectra'). 
As we saw for $K$-theory,  it is then especially easy to calculate $[\Sigma^{\infty}X, E]$ since
$$
[\Sigma^{k+1}X,E_{k+1}]_{\un}\cong [\Sigma^kX,\Omega E_{k+1}]_{\un}
\cong [\Sigma^kX,E_k]_{\un}
$$
and all maps in the limit system are isomorphisms.

In particular
$$
\pi_n(E)=\pi_n(E_0) \quad \text{ for } n\ge0
$$
and in fact more generally
$$
\pi_n(E)=\pi_{n+k}(E_k)\quad \text{ for } n+k\ge0.
$$

(b) If $X$ is a $\Omega$-spectrum,  the 0th term $X_0$ has the remarkable property
that it is equivalent to a $k$-fold loop space for each $k$ (indeed, $X_0 \simeq \Omega^kX_k$). 
Spaces with this property are called {\em $\Omega^{\infty}$-spaces} 
(sometimes pronounced `infinite loop spaces'). The space $X_0$ does not retain information 
about negative homotopy groups of $X$, but if $\pi_i(X)=0$ for $i<0$ 
(we say $X$ is {\em connective}), and  we retain information about {\em how} it is a $k$-fold 
loop space for each $k$ we have essentially recovered the spectrum $X$. 
The study of  $\Omega^{\infty}$-spaces is equivalent to the category of 
connective spectra in a certain precise sense.
\end{remark}

To get the best formal behaviour,  we impose an even stronger condition than being a 
$\Omega$-spectrum.
 
\begin{defn}
\label{defn:Mayspectrum}
A {\it May spectrum\/} is a spectrum so that the adjoint bonding maps 
$$
\tilde \sigma\colon X_k\stackrel{\cong } \lra  \Omega X_{k+1}
$$
are all homeomorphisms.

\end{defn}

\begin{remark}
Spectra in this strong sense are rather rare in nature, but there is a left adjoint
$$
L\colon \text{Spectra} \to \text{May spectra} 
$$
to the inclusion of May spectra in spectra. On reasonable spectra 
(including those for which the bonding maps are cofibrations) it is given by
$$
(LE)_k=\colim_{ s} \Omega^sE_{k+s}.
$$
For instance
$$
(L\Sigma^\infty X)_k=\colim_{ s}\Omega^s\Sigma^{k+s} X
$$
and we have a version of the $Q$-construction $(L\Sigma^\infty X)_0=QX$.

In general we will omit mention of the functor $L$, for example writing
$\Sigma^{\infty}X$ for the spectrum associated to the suspension 
spectrum and $\bbS$ for the 0-sphere May spectrum.
\end{remark}

One can then proceed with homotopy theory of May spectra very much as with 
spaces or forming  the derived category. One wants to invert 
$\pi_*$-isomorphisms and work with
$$
\text{Spectra}[(\pi_*\text{-isos})^{-1}].
$$
To avoid set-theoretic difficulties with categories of fractions, we construct
this homotopy category directly. First we define cells and spheres using shifted 
suspension spectra and then CW-spectra.
Since cells are compact in a suitable sense,  it is elementary to form
 CW-approximations. For any spectrum $E$ we may construct
a spectrum $\Gamma E$ from cells, together with a map 
$$
\Gamma E \to E
$$
which is  a weak equivalence. It is then a formality that $\Gamma $
provides a functor in the homotopy category, and it is called the 
{\em CW-approximation functor}. Using this construction, we find
$$
\text{Spectra}[(\pi_*\text{-isos})^{-1}]\simeq \text{Ho}(\text{CW-Spectra})
$$
and this is usually just called the {\em homotopy category of spectra}.

This has good formal  properties like the derived category. It is triangulated,
has products, sums and  internal homs (function spectra).

\section{The smash product (Step 2).}
\label{sec:smash}

We have now completed Step 1 by constructing a suitable additive 
category, and we now proceed to Step 2 and  endow the category of spectra
with additional structure, 
especially that of an associative and commutative smash product. 
This is made a little harder because it is necessary to restrict or 
otherwise adapt the category of spectra that we have found so far. 

We would like to form a smash product  $E\wedge F$ of spectra $E$ and $F$
from the terms $E_k\wedge F_l$ in some way. In the first instance, 
we have a doubly indexed collection of spaces, and to make a spectrum 
out of it we would need to somehow combine all possibilities or select 
from them.
If done too  naively, we lose all hope of associativity of the result. 
There are several approaches to avoiding this problem. We describe three: 
the EKMM approach, the approach via symmetric spectra, and that via orthogonal 
spectra. We emphasize that these
all give derived categories which are equivalent in a very strong sense 
\cite{MMSS}, but as usual each has its
own advantages and disadvantages. In each case there is a sphere spectrum
$\bbS$ which is a ring (using the smash product) and the spectra are
modules over $\bbS$.

We begin with the EKMM approach for the same
reason one starts homotopy theory with spaces rather than simplicial sets, 
but (partly because of what is omitted in this account) I suspect that 
commutative algebraists will prefer the symmetric spectra described in 
Subsection \ref{subsec:symmetric} below.

\subsection{Method 1: EKMM spectra.} The acronym refers to 
Elmendorf, Kriz, Mandell and May \cite{EKMM}. They call their category 
of spectra $\bbS$-modules, where $\bbS$ is the sphere spectrum, but 
this name also describes other categories, so we refer to `EKMM spectra'.

\vskip .1in

First,  there is a partial solution based on not making choices, sometimes
called {\em coordinate free spectra}. We extend
the notation for spheres and suspensions to permit arbitrary real vector
spaces, so that $S^V$ denotes the one-point compactification of $V$
and $\Sigma^VX:=X \sm S^V$.

\begin{defn}
\label{defn:spectraspectra}
(i) A {\it universe\/} is a countable dimensional real inner product space. 
An {\em indexing space} in a universe $\uc$
is a finite dimensional  sub inner product space 
$V\subseteq \uc$.

(ii) A {\it spectrum $E$ indexed on $\uc$\/} is a collection
of spaces $E_V$ where $V$ runs through indexing spaces $V$ in 
$\uc$  together with a transitive system of bonding maps
$$
\sigma_{V,W}:  \Sigma^{W-V} E_V \to E_W
$$
whenever $V\subseteq W$, where $W-V$ denotes the orthogonal complement of 
$V$ in $W$.

(iii) Such a spectrum is a {\it May spectrum\/} if all adjoint bonding maps 
 $\tilde\sigma \colon E_V \overset{\cong}{\to} \Omega^{W-V}E_W$  are 
homeomorphisms.
\end{defn}

\begin{remark}
(a) From any cofinal sequence of indexing spaces one may fill in gaps by using 
suspensions. Hence we consider a spectrum to be specified by such a 
cofinal sequence.

For example, if we choose a cofinal sequence
$\re\subseteq \re^2\subseteq\cdots\subseteq \uc$ with $n$ corresponding
to $\R^n$ we can convert a spectrum as in \ref{defn:spectra}
 into a spectrum indexed on $\uc$. 

(b) We may also change universes. If $f: \uc \lra \vc$ is an isometry, 
we may use $f$ to convert a spectrum $E$ indexed
on $\uc$ to a spectrum $f_*E$ indexed on $\vc$, by taking $(f_*E)(V):=E(f^{-1}V)$.
\qqed
\end{remark}

\begin{defn}
Given a spectrum  $E$ indexed on $\cU$ and a spectrum $F $ indexed on $\cV$, 
one may define the external smash product $E\extsmash F$ indexed on 
$\uc\oplus\vc$ by taking
$$
(E \extsmash F)(U\oplus V) := E(U)\wedge F(V)
$$
on the cofinal sequence of indexing spaces of the form $U\oplus V$.
\end{defn}

The merit of the definition is that no choices are involved. 
Thus if $G$ is a spectrum indexed on $\wc$,  
there is a coherent natural associativity isomorphism
$$
(E\ov{\wedge} F) \ov\wedge G\cong E\ov\wedge(F\ov\wedge G)
$$
of spectra indexed on $\uc\oplus \vc\oplus \wc$.

The problem is that if $E$ and $E'$ are both indexed on $\uc$ then 
$E\ov\wedge E'$ is indexed on $\uc \oplus \uc$ rather than on $\uc$ itself.
The old fashioned solution is to {\it choose\/} an isometric isomorphism
$$
i\colon\uc\oplus \uc \overset{\cong}\lra  \uc, 
$$
and use it to index $E \extsmash E'$ on $\uc$: we define
$$
E\wedge_iE':= i_*(E\ov\wedge E').
$$
This depends on $i$, but because the space $\lc (2) :=\lc(\uc\oplus\uc,\uc)$
of linear isometries is contractible,  the choice of $i$ is relatively 
unimportant,  and because the spaces $\lc (n):=\lc (\uc^{\oplus n}, \uc)$
are contractible for $n \geq 1$, this gives a coherently commutative and associative 
operation {\em in the homotopy category}. This method of internalizing the smash 
product is quite useful, but to obtain the good properties before passing to homotopy
we must work a little harder.

The EKMM solution is to use {\it all\/} choices. The key to this is the 
twisted half-smash product construction, which we only describe in general terms.

\begin{construction}
Given 

(i) a space $A$, 

(ii) a map $\alpha : A{\to} \lc(\uc,\vc)$, and

(iii) a spectrum $E$ indexed on $\uc$, 

\noindent
we may form the {\em twisted half-smash product} $A\ths E$. This is
a spectrum indexed on $\vc$ formed by  assembling the 
spectra $\alpha (a)_*E$ for all $a \in A$.

The twisted half-smash product is natural for maps of $A$ and $E$.
It is also homotopy invariant in the strong sense that 
the homotopies need not be compatible with the structure maps $\alpha$.
\end{construction}

\begin{example}
(a) If we choose the one point space, we recover the earlier change
of universe construction. If 
 $A=\{i\}\subseteq \lc(\uc,\vc)$ then $\{i\}\ltimes E=i_* E$

\vskip .1in

\noindent
(b) If we take $A=\lc(\uc^{\oplus 2}, \uc)$ and let $\alpha$ 
 be the identity we obtain a canonical way to internalize a 
smash product. We may take 
$$
E\wedge' E':= \lc(\uc^{\oplus 2}, \uc)\ltimes (E\ov\wedge E').
$$
By the naturality, all choices of  $\wedge_i$ are contained in this, but it 
is still a bit too big to be associative. \qqed
\end{example}

Restricting attention to spectra with a little extra structure, 
one may remove some flab from this smash product and make 
an associative one.

\begin{defn}
An $\lb$-spectrum is a May spectrum $E$ with an action $\lb E \lra E$, 
where $\lb$ is  the functor defined by 
$\lb E:=\lc(1) \ltimes E$. We may view this as a continuous family of maps 
$f_*E\to E$ where $f\in \lc(\uc,\uc)$,  compatible with composition.
\end{defn}

There are plenty of examples of $\lb$-spectra. For example the sphere spectrum
$\bbS$ is an $\lb$-spectrum, as is any suspension spectrum. In general, 
any spectrum $E$, is homotopy equivalent to the $\lb$-spectrum $\lb E$ 
(since $\lc (1)$ is contractible).

\begin{defn}
The smash product of $\lb$-spectra $M$, $N$ is then defined by 
$$
M\wedge_{\lc}N:= \lc(2)\ltimes_{\lc(1)\times\lc(1)} (M\smbar N).
$$
More precisely, it is the coequalizer
$$
\diagram
(\lc(2)\times\lc(1)\times\lc(1)) \ltimes(M\smbar N) \rto<0.4ex> \rto<-0.4ex>
& \lc(2)\ltimes(M\smbar N) \rto 
&M\wedge_{\lc} N
\enddiagram
$$
using the maps
$$
(\theta,\varphi,\psi)
 \longmapsto (\theta\circ(\varphi\oplus \psi), M,N) 
$$
and 
$$
(\theta,\varphi,\psi)
 \longmapsto (\theta, \varphi_*M,\psi_* N).
$$
\end{defn}

This finally  gives a good smash product. 

\begin{prop}
(Hopkins) The smash product $\sm_{\lc}$ is commutative and associative.
\end{prop}

\begin{remark} 
\label{iteratedsmash}
This proposition is a formal consequence of two key features of $\lc$:

(1) \, $\lc(i+j)\cong \lc(2) \times_{\lc(1)\times\lc(1)} \lc(i)\times \lc(j)$

(2) \, $\lc(2)/\lc(1)\times\lc(1)=*$.

Building on these, we may also rearrange the iterated product
$$
M_1\wedge_{\lc}\cdots \wedge_{\lc} M_n \cong\lc(j) \ltimes_{\lc(1)^j} 
(M_1\smbar \cdots\smbar M_n).
$$
This is useful in recognizing monoids and commutative monoids.
\end{remark}

It is convenient to ensure that $\bbS$ is itself the unit for the smash 
product, so we restrict attention to the category of $\bbS$-modules
(i.e., $\lb$-spectra for which the natural weak equivalence
$\bbS\wedge_{\lc}M\overset{\simeq}{\to} M$ is actually an isomorphism). Since
every $\lb$-spectrum $E$ is weakly equivalent to the $\bbS$-module
$\bbS \sm_{\lc} E$, and since the smash product preserves $\bbS$-modules, 
this is no real restriction.







\subsection{Method 2: symmetric spectra.}
\label{subsec:symmetric}

This method is due to Jeff Smith, with full homotopical details  published in 
\cite{HSS}. It gives a more elementary and combinatorial
construction of a symmetric monoidal category of spectra, 
but the construction of the {\em homotopy} category is much more indirect 
and requires fluency with Quillen model categories.
This is directly analogous to the situation for spaces. Most people find
it more intuitive to work with actual topological spaces with homotopies
being continuous one-parameter  families of maps, and to restrict to  CW-complexes
to obtain a well-behaved homotopy category. However one may construct
the homotopy category using  simplicial sets instead. This gives a purely 
combinatorial model with some superior formal properties, but the construction
of the homotopy category requires considerable work. Because of these superior
properties, it is usual to base symmetric spectra on simplicial sets 
(i.e.,  in Step 0) rather than on topological spaces.

\begin{defn}
(a)  A {\it symmetric sequence\/} is a sequence
$$
E_0,E_1,E_2,\dots,
$$
of pointed simplicial sets with basepoint
 preserving action of the symmetric group $\Sigma_n$ on $E_n$.

\vskip .1in

\noindent
(b) We may define a tensor product $E \tensor F$ of symmetric sequences $E$ and
$F$ by  
$$
(E\otimes F)_n:= \bigvee_{p+q=n} (\Sigma_{n})_+ \wedge_{\Sigma_p\times \Sigma_q} (X_p\wedge Y_q), 
$$
where the subscript $+$ indicates the addition of a disjoint basepoint.
\end{defn}

It is elementary to check that this has the required formal behaviour.

\begin{lemma}
The product $\otimes$ is symmetric monoidal with unit
$$
(S^0,*,*,*,\dots).
$$
\end{lemma}

\begin{example} The sphere is the symmetric sequence 
$\bbS:=(S^0,S^1,S^2,\dots)$. Here  $S^1=\Delta^1/\po \Delta^1 $ is the simplicial circle 
and the higher simplicial spheres are defined by taking smash powers, 
so that  $S^n=(S^1)^{\wedge n};$ this also explains the actions of the symmetric groups.

It is easy to check that {\em the sphere  is a commutative monoid}  in the category of 
 symmetric sequences.
\end{example}

\begin{defn}
A {\it symmetric spectrum $E$\/} is a left $\bbS$-module in symmetric sequences.

Unwrapping the definition, we see that this means $E$ is given by 

(1) a sequence $E_0,E_1,E_2,\dots$ of simplicial sets, 

(2) maps $\sigma\colon S^1\wedge X_n\to X_{n+1}$, and 

(3) basepoint preserving  left actions of $\Sigma_n$ on $X_n$ which are 
compatible with the actions in the sense that the composite maps 
$S^p\wedge X_n \to X_{n+p}$ are $\Sigma_p \times \Sigma_n$ equivariant.
\end{defn}

\begin{defn}
The smash product of symmetric spectra
is
$$
\diagram
E \sm_{\bbS} F:=\coeq (E \otimes \bbS \otimes F \rto<0.4ex> \rto<-0.4ex>& E \otimes F).
\enddiagram$$
\end{defn}

\begin{prop}
The tensor product over $\bbS$ is a symmetric monoidal product on the category of symmetric
spectra.
\end{prop}

It is now easy to give the one example most important to commutative 
algebraists.
\begin{example}
For any abelian group $M$, we define the Eilenberg-MacLane symmetric spectra.
For a set  $T$ we write $M \otimes T$ for the $T$-indexed sum of copies
of $M$; this is natural for maps of sets and therefore extends to an 
operation on simplicial sets. We may then define the Eilenberg-MacLane
symmetric spectrum
$HM:=(M\otimes S^0,M\otimes S^1, M\otimes S^2,\dots)$. 
It is elementary to check that if $R$ is a commutative ring, then  $HR$ is a monoid 
in the category of $\bbS$-modules, and if $M$ is an $R$-module, $HM$ is a module over
$HR$.
\end{example}

We will not spoil the impression of immediate accessibility of symmetric spectra by
explaining how to form the associated homotopy category: 
one needs to restrict to a good class of symmetric spectra and then invert a certain 
collection of weak equivalences. The weak equivalences are not just homotopy 
isomorphisms, so  this involves some work in the framework of model categories.

\subsection{Method 3: orthogonal spectra.} Combining the merits of EKMM spectra and
symmetric spectra there is a third option \cite{MMSS, MM}.

For this we let $\cI$ denote the category of finite dimensional real inner product
spaces; the set of morphisms between a pair of objects forms a topological space, 
and the composition maps are continuous. For example $\cI (U,U)$ is the orthogonal 
group $O(U)$.

\begin{defn}
 An {\em $\cI$-space} is a continuous functor $X: \cI \lra \basedspaces$ to the category 
of based spaces.
\end{defn}

Notice  the large amount of naturality we require: for example $O(U)$ acts on 
$X(U)$, and an isometry $U \lra V$ gives a splitting $V=U \oplus V'$ so that
$X(U) \lra X(V)= X(U\oplus V')$ is also $O(U)$-equivariant. 

 A very important  example is the functor $\bbS$ which takes an inner product
space $V$ to its one point compactification $S^V$.

There is a natural external smash product of $\cI$-spaces, so that if $X$ and 
$Y$ are $\cI$-spaces we may form
$$
X \smbar Y : \cI \times \cI \lra \basedspaces
$$
by taking $(X \smbar Y) (U,V):= X(U) \sm Y(V)$.

\begin{defn}
An {\em orthogonal spectrum} is an $\cI$-space $X$ together with a natural map 
$$
\sigma : X \smbar \bbS \lra X \circ \oplus
$$
so that the evident unit and associativity diagrams commute. Decoding this, 
we see that the basic structure consists of maps
$$
\sigma_{U,V}: X(U) \sm S^V \lra X(U \oplus V), 
$$
and this commutes with the action of $O(U) \times O(V)$.
\end{defn}


One may define the  objects which play the role of rings without defining
the smash product. 
 
\begin{defn} 
An {\em $\cI$-functor with smash product} (or $\cI$-FSP) is an $\cI$-space $X$ with 
a unit $\eta : \bbS \lra X$ and a natural map $\mu : X \smbar X \lra X \circ \oplus$. 
We require that $\mu$ is associative, that $\eta $ is a unit (and central) in the
evident sense. For a commutative $\cI$-FSP we impose a commutativity condition on 
$\mu$.
\end{defn}

Note that the unit is given by maps
$$
\eta_V : S^V \lra X(V)
$$
and the product $\mu$ is given by  maps
$$
\mu_{U,V}: X(U) \sm X(V) \lra X(U\oplus V).
$$
Thus, by composition we obtain maps
$$
X(U) \sm S^V \lra X(U) \sm X(V) \lra X(U \oplus V), 
$$
and one may check that these give an $\cI$-FSP the structure of an 
orthogonal spectrum.

\begin{remark}
The notion of $\cI$-FSP is closely related to the FSPs introduced by  
B\"okstedt in algebraic $K$-theory before a symmetric monoidal smash 
product was available. An FSP is a  functor from simplicial sets to 
simplicial sets with unit and product. The restriction 
of an FSP to (simplicial) spheres is analogous to a $\cI$-FSP and gives 
rise to a ring in symmetric spectra. \qqed
\end{remark}

To define a smash product one first defines the smash product of 
$\cI$-spaces by using a Kan extension to internalize the 
 product $\smbar$ described above. Now observe that $\bbS$ is a monoid for 
this product and define the smash product of orthogonal spectra to be the 
coequalizer
$$
\diagram
X \sm_{\bbS} Y:=\coeq (X \sm \bbS \sm Y \rto<0.4ex> \rto<-0.4ex>& X \sm Y).
\enddiagram$$
The monoids for this product are essentially the same as $\cI$-FSPs. 

As for symmetric spectra, a fair amount of model categorical work is
necessary to construct the associated homotopy category, but orthogonal 
spectra have the advantage that the weak equivalences are the homotopy 
isomorphisms.

\section{Brave new rings.}

Once we have a symmetric monoidal product on our chosen category of spectra 
we can implement the dream of the introduction: choose  a ring spectrum $\boldR$
(i.e., a monoid in the category of spectra), 
form the category of $\boldR$-modules or $\boldR$-algebras and 
then pass to homotopy. We may then attempt to use algebraic methods and
intuitions to study $\boldR$ and its modules. We use bold face for ring
spectra to remind the reader that although the methods are familiar, we are
not working in a conventional algebraic context. The `brave new ring' terminology 
is due to Waldhausen, and nicely captures both the wonderful possibilities and
the denaturing effect of inappropriate generality. Some in  the new wave
prefer the term `spectral ring'.

In turning to examples, we remind the reader that the equivalence results of
\cite{MMSS} mean that we are free to choose the category most convenient for
each particular application. 

\begin{example}
If we are prepared to use symmetric spectra,  we already have the example of
the Eilenberg-MacLane spectrum $\boldR=HR$ for a classical commutative ring $R$.

The construction of the Eilenberg-MacLane symmetric spectra  gives 
a  functor $R\text{-modules } \longrightarrow HR\text{-modules}$ and 
passage to homotopy groups gives a functor Ho($HR$-mod) $\to R$-modules.
It is much less clear that there are similar comparisons of derived categories but 
in fact the derived categories are equivalent.

\vskip .1in

\begin{thm} (Shipley \cite{ShipleyQuillen}) There is a Quillen equivalence
between the category of $R$-modules and the category of $HR$-modules,
and hence in particular a triangulated equivalence
$$
D(R)=\text{Ho}(R\text{-modules}) \simeq \text{Ho}(HR\text{-modules})=D(HR)
$$
of derived categories. 
More generally, one may associate a ring spectrum $HR$ to any DG ring
$R$, so that $H_*(R)=\pi_*(HR)$,  and the same result holds. \qqed
\end{thm}

\end{example}

Thus working with spectra does recover the classical algebraic
derived category. However there are plenty more examples.

\begin{example} 
\label{eg:cochains}
For any space $X$ and a commutative ring $k$ we may form the function spectrum 
$\boldR=\text{map}(\Sigma^\infty X, Hk)$. It is obviously  an $Hk$-module, 
but using the diagonal on $X$ it is also a commutative $Hk$-algebra.
Certainly 
$$
\pi_*(\text{map}(\Sigma^{\infty}X,Hk))=H^*(X;k), 
$$
and $\boldR$ should be viewed as a commutative substitute for the DG algebra of 
cochains $C^*(X;k)$. Similarly, a map  $Y \lra X$ makes the substitute for 
$C^*(Y;k)$ into an $\boldR$-module. 
The commutative algebra of this ring spectrum $\boldR$ is extremely interesting 
(\cite{DGI}) and discussed briefly in Section \ref{sec:localring}.
\end{example}

\begin{example}
If $G$ is a group or a monoid. Then
$$
\boldR=\Sigma^\infty G_+
$$
is a monoid, commutative if $G$ is abelian. The case  $G=\Omega X$ for a space $X$ 
is important in geometric topology (here one should use Moore loops to ensure that
$G$ is strictly associative).
\end{example}

\begin{example}
We may apply  the algebraic $K$-theory functor to any ring spectrum $\boldR$ to 
form a spectrum $K(\boldR)$. If $\boldR$ is a commutative ring 
spectrum so is $K(\boldR)$.

This generalizes the classical case in the sense that $K(HR)=K(R)$ (where 
the right hand side is the version of algebraic $K$-theory based on finitely 
generated free modules). Another important example comes from 
geometric topology: $K(\Sigma^{\infty}\Omega X_+)$ is Waldhausen's $A(X)$ \cite[VI.8.2]{EKMM}.
The spectrum $A(X)$ embodies
 a fundamental step in the classification of manifolds
\cite{Waldhausen}.
The calculation of its homotopy groups can often be approached using the 
methods described for algebraic $K$-theory in Subsection \ref{subsec:algK}.
\end{example}

To import many of the classical examples we need to decode what is needed
to make a commutative $\bbS$-algebra in the EKMM sense, 
 using  Remark \ref{iteratedsmash}. 

\begin{lemma}
\cite[II.3.6]{EKMM}
A commutative $\bbS$-algebra is essentially the same as an $E_\infty$-ring 
spectrum i.e., a spectrum $X$ with maps
$$
\lc(\uc^k,\uc) \wedge X^{\wedge k} \to X
$$
with suitable compatibility properties. More precisely, if $X$ is an $E_{\infty}$-ring
spectrum, the weakly equivalent EKMM-spectrum $\bbS \sm_{\lc}X$ is a commutative
$\bbS$-algebra.
\end{lemma}

\begin{remark}
(i) The space $\lc(\uc^k,\uc)$ is $\Sigma_k$-free and contractible, and taken together
these spaces form the {\em linear isometries operad}. Any other sequence
$$
\oc(0),\oc(1),\oc(2),\dots
$$
of contractible spaces with free actions of symmetric groups and similar 
compositions is called an 
$E_\infty$-operad \cite{Mayoperads}. Up to suitable equivalence,  it does not depend 
which $E_\infty$-operad is used, so that although the linear isometries 
operad is rather special because of \ref{iteratedsmash}, using it results in 
no real loss of generality.

(ii) This method allows an obstruction theoretic approach to constructing
$\bbS$-algebra structures, where the obstruction groups are based on a
topological version of Hochschild cohomology (or a topological version 
of Andr\'e-Quillen cohomology in the commutative case).
\end{remark}

\begin{cor} The following spectra are commutative $\bbS$-algebras:
the bordism spectra $MO$ and  $MU$, the $K$-theory spectrum $K$  
and its connective cover $ku$.
\end{cor}

\vskip .2in

\noindent
{\bf Proof for $MO$:} We may use the Grassmann model for the classifying
space $BO(N)$. In fact for a universe $\uc$ we may take
$BO(N)=Gr_{N}(\uc), $ the space of $N$-dimensional subspaces of 
$\uc$. Noting that $\uc \cong U \oplus \uc$ for any indexing subspace $U$, 
we have natural maps
$$
\begin{array}{ccc}
MO(U)_{\uc} \wedge MO(V)_{\uc} & & MO(U\oplus V)_{\uc\oplus \uc} \\
||  & &  || \\
Gr_{|U|}(U\oplus \uc)^{\ga_{|U|}} \wedge Gr_{|V|} (V\oplus \uc)^{\ga_{|V|}} & \longrightarrow & Gr_{|U\oplus V|} (U\oplus V\oplus \uc \oplus\uc)^{\ga_{|U\oplus V|}}
\end{array}
$$

A choice of isometry $\uc^{\oplus 2} \to \uc$ gives a map 
 $MO(U\oplus V)_{\uc \oplus \uc} \lra MO(U\oplus V)_{\uc}$, 
and assembling these we obtain a map
$$
\lc(\uc^{\oplus 2},\uc) \wedge MO_{\uc} \wedge MO_{\uc} \to MO_{\uc}
$$
and similarly for other numbers of factors. 

Another way to construct $MO$ as  a commuative $\bbS$-algebra is as
an $\cI$-FSP. Indeed we may take $MO'(V)$ to be the Thom space of 
the tautological bundle over $Gr_{|V|}(V \oplus V)$, and then the 
structure maps are constructed just as above. The inclusions 
$$Gr_{|V|}(V \oplus V)^{\gamma_V} \lra Gr_{|V|}(V \oplus \uc)^{\gamma_V}$$
give rise to a map $MO' \lra MO$ of the associated spectra. Since the
maps of spaces 
become more and more highly connected as the dimension of $V$ increases, 
 this shows that $MO' \simeq MO$.
\vskip .2in

\noindent
{\bf Conclusion:} There are many examples of commutative $\bbS$-algebras.

\section{Some algebraic uses of ring spectra.}

The main purpose of this article is to introduce spectra, but we want
to end by showing they are useful in algebra.  
Our principal example of commutative algebra is in the next section, 
but we mention a number of other applications briefly here.

\subsection{Topological Hochschild homology and cohomology.}

Given a $k$-algebra $R$ with $R$ flat over $k$, 
we may define the  Hochshild homology and cohomology using homological algebra over 
$R^e:=R \tensor_k R^{op}$, by taking
$$
HH_*(R|k):=\mathrm{Tor}^{R^e}_*(R,R)
$$ 
and 
$$
HH^*(R|k):=\mathrm{Ext}_{R^e}^*(R,R);
$$
we have included $k$ in the notation for emphasis, but it is often 
omitted. We may make precisely parallel definitions for ring spectra. 
In doing so, we emphasize that all Homs  of ring spectra in this article
are derived Homs (sometimes written $\boldR\Hom$) and all tensors of
ring spectra are derived (sometimes written $\tensor^{\bf L}$). 
Because of this, it is no longer necessary to make a flatness hypothesis. 
If $\boldR$ is
a $\boldk$-algebra spectrum we may  define the topological versions
using homological algebra over the ring spectrum 
$\boldR^e:=\boldR \sm_{\boldk}\boldR^{op}$, defining
the   Hochschild homology spectrum by
$$THH_{\bullet} (\boldR | \boldk):=\boldR \sm_{\boldR^e}\boldR$$
and the topological Hochschild cohomology spectrum by
$$THH^{\bullet} (\boldR |\boldk ):=\Hom_{\boldR^e}(\boldR,\boldR).$$
The $\bullet$ subscript and superscript indicates whether homology
or cohomology is intended. When $\boldk$ is omitted in the notation for 
$THH$, it is assumed to be the sphere spectrum $\boldk =\bbS$; in this 
case $THH$ was first defined by B\"okstedt by other means before
good categories of spectra were available. 
We may obtain purely algebraic topological Hochschild 
homology and cohomology groups by taking homotopy, so that 
$THH_*(\boldR | \boldk )=\pi_*(THH_{\bullet}(\boldR | \boldk))$
and $THH^*(\boldR | \boldk )=\pi_*(THH^{\bullet}(\boldR | \boldk))$.
Alternative notations such as $THH (\boldR | \boldk )=
THH_{\bullet}(\boldR | \boldk )= THH^{\boldk}(\boldR)$
and $THC (\boldR | \boldk )=THH^{\bullet}(\boldR | \boldk)
=THH_{\boldk}(\boldR)$ also occur in 
the literature, but unfortunately $THC$ may be confused with cyclic homology.

Under flatness hypotheses to ensure $\pi_*(\boldR^e)=(\pi_*(\boldR))^e$, 
there are spectral sequences 
$$HH_*(\pi_*(\boldR)| \pi_*(\boldk) ) \Rightarrow \pi_*(THH_{\bullet} (\boldR | \boldk))$$
and
$$HH^*(\pi_*(\boldR) | \pi_*(\boldk)) \Rightarrow \pi_*(THH^{\bullet} (\boldR | \boldk)).$$
In particular if $\boldR=HR$ and $\boldk=Hk$ for a conventional rings $R$ and
$k$ with $R$ flat over $k$,  the spectral sequences collapse for dimensional reasons to show
that the Hochschild homology and cohomology of $R$ is equal to the 
topological Hochschild homology and cohomology of $HR$. 

Two uses of the Hochschild groups are to provide invariants for algebraic
$K$-theory and to provide an obstruction theory for extensions of
rings; both of these applications have parallel versions in the topological 
theory. We briefly describe some applications below. There is also a topological version of
Andr\'e-Quillen cohomology \cite{RW,Basterra,BasterraMcCarthy}
which can be used to give an obstruction theory for extensions
of commutative ring spectra.

\subsection{Algebraic $K$-theory and traces.}
\label{subsec:algK}
The algebraic $K$-theory $K_*(R)$ of a ring $R$ is notoriously hard to 
calculate, and one method is to use 
trace maps to attempt to detect $K$-theory. B\"okstedt, Hesselholt, Madsen and others have
calculated the  $p$-complete algebraic $K$-theory of suitable $p$-adic rings 
\cite{BM,HM,Madsen} using spectral refinements of classical traces. 
The relevant constructions were first made  using B\"okstedt's FSPs.

The classical Dennis trace map $K_*(R) \lra HH_*(R)$ lands in the Hochschild 
homology of $R$, and B\"okstedt
has given a topological version, which is a map  $K(R) \lra THH_{\bullet} (HR|\bbS )$ of 
spectra. Taking homotopy of B\"okstedt's map gives a refinement of the Dennis trace.
However there is more structure to exploit: the cyclic structure of the 
Hochschild complex gives 
a circle action on $THH_{\bullet}(HR | \bbS )$ and the geometry of B\"okstedt's map 
shows it has equivariance properties. The fixed point spectra
 $THH_{\bullet}(HR | \bbS)^C$ for finite cyclic groups $C$ are related in the
usual way, but also by maps arising from the special `cyclotomic'
nature of the Hochschild complex; taking both structures into account, one may
construct a {\em topological cyclic homology spectrum} $TC(R)$ from these
fixed point spectra. The 
construction of  $TC (R)$ from $THH_{\bullet}(HR | \bbS )$ can be modelled 
algebraically, and this makes the homotopy groups $TC_*(R)$ 
relatively accessible to calculation. Because the relationships between 
fixed point sets correspond to structures in algebraic $K$-theory, 
B\"okstedt, Hsiang and Madsen \cite{BHM} are able to construct a map
$$trc: K(R) \lra TC(R) $$
of spectra. Again we may take homotopy to give the {\em cyclotomic trace} 
$K_*(R) \lra TC_*(R)$. 
This is a very strong invariant, and for certain classes of rings $R$ it is 
actually a $p$-adic isomorphism. Indeed, McCarthy \cite{McCarthy} has shown
that the cyclotomic trace always induces a profinite isomorphism of 
{\em relative} $K$-theory. From the known $p$-adic algebraic $K$-theory of 
perfect fields $k$ of characteristic $p>0$, Madsen and Hesselholt 
\cite{HM} deduce that the cyclotomic trace is a $p$-adic isomorphism in 
degrees  $\geq 0$ whenever $R$ is an algebra over the Witt vectors $W(k)$ 
which is finite as a  module. This, combined with calculations of $TC_*(R)$
has been used to calculate $K_*(R)_p^{\wedge}$ for many complete 
local rings $R$, including $R=\Z_p^{\wedge}$ and truncated polynomial rings
$ k[x]/(x^n)$, and to prove  the Beilinson-Lichtenbaum
conjectures  on the $K$-theory of Henselian discrete valuation fields
of mixed characteristic \cite{Hesselholt}.

\subsection{Topological equivalence.}

Two rings are said to be {\em derived equivalent} if their derived categories
are equivalent as triangulated categories. The best known example is that of
Morita equivalence, showing that a ring is derived equivalent to the ring
of $n \times n$ matrices over it. Since useful invariants can be constructed from 
the derived category, the freedom to replace a ring by a derived equivalent 
ring can be very useful.

For ring spectra, it is natural to consider also the stronger condition that
the module categories are Quillen equivalent (this implies that
their derived categories are triangulated equivalent, but it is usually 
a stronger condition). We then say that the ring spectra are 
{\em Quillen  equivalent}.

Just as any derived equivalence of rings is given by tensoring with a complex
of bimodules, any Quillen  equivalence between 
 ring spectra is given by smashing with a bimodule spectrum 
\cite{SchwedeShipley}. In particular, any Quillen  equivalence between 
DG algebras is given by smashing with a bimodule spectrum, but  Dugger and Shipley 
\cite{DuggerShipley} have given an example to show that it need not be given
by tensoring with a complex of bimodules. Based on work of Schlichting, 
they have also given an example to show that derived equivalent ring spectra
need not be Quillen  equivalent (although derived equivalence and 
Quillen  equivalence agree for ungraded rings).

Two DG algebras are {\em quasi-isomorphic} if  they are related by a chain of homology 
isomorphisms. Similarly, two ring spectra are {\em topologically equivalent} if 
they are related by a chain of homotopy isomorphisms. If the DG algebras 
are quasi-isomorphic, 
the associated ring spectra are topologically isomorphic, but Dugger and Shipley
have given an example to show that topological equivalence does not imply 
quasi-isomorphism. Perhaps the best way to think about this is 
 that there is a ring map $\bbS \lra H\Z$; viewing a DG $\Z$-algebra
as an $H\Z$-algebra, we may view  it  as a $\bbS$-algebra by 
restriction. It is then not surprising that an equivalence of $\bbS$-algebras need
not be an equivalence of $H\Z$-algebras. 
Since topological equivalence implies Quillen  equivalence, 
this shows that viewing DG algebras as ring spectra can have useful consequences. 

There is an obstruction theory for extensions of rings based on Hochschild cohomology, 
and a parallel theory for extensions of ring spectra  based on topological 
Hochschild cohomology. The Dugger-Shipley example is based on the comparison between
algebraic and topological Hochschild cohomology.

\section{Local ring spectra.}
\label{sec:localring}

Finally we turn to  the spectral analogue of a commutative Noetherian
local ring $R$ with residue field $k$.
In effect we are extending the idea of trying to do commutative algebra 
entirely in the 
derived category. When notions can be reformulated in these terms, we 
gain considerable flexibility.

We consider a map $\boldR \lra \boldk$ of commutative ring spectra, viewed as an 
analogue of the map  from a commutative local ring $R$ to its residue field $k$.
One example is to take $\boldR=HR \lra Hk=\boldk$, and we refer to this as the
local algebra example.  A second example is to take 
$\boldR=C^*(X;k)$ (in the sense of Example \ref{eg:cochains}) for a space 
$X$ and $\boldk=Hk$, and we refer to this as the example of cochains on a space. 

We study the map $\boldR \lra \boldk$ with the eyes of commutative algebra, and
illustrate the possibilities with results from \cite{DGI}. The projects
of Waldhausen \cite{Waldhausenchrom} and Rognes \cite{Rognes} to give an 
analysis of chromatic stable homotopy theory by using commutative algebra 
and Galois theory are beyond the scope of these notes.

\subsection{Regularity.}
Serre's characterization of regularity states that if $R$ is a commutative
Noetherian local ring with residue field $k$,  then $R$ is regular if and only if 
$k$ has a finite free resolution by $R$-modules. In the derived category 
we can construct the resolution as a complex of $R$-modules in finitely 
many steps from $R$ by completing triangles and passing to direct summands
(we say ``$k$ is finitely built from $R$''). This leads to our definition. 

\begin{defn}
We say that $\boldR \lra \boldk$ is {\em regular} if $\boldk$ can be finitely built from $\boldR$.
\end{defn}

Serre's characterization shows that $R$ is regular  if and only if 
$\boldR=HR \lra Hk=\boldk$ is regular as a ring spectrum, so in the commutative
algebra example the new notion agrees  with the classical one.

Regularity is an interesting condition for many other examples.
In the case of cochains on a space, we consider the special case 
with $k=\bF_p$ for some prime $p$, and $X$ $p$-complete. Thus  
 $\boldR=C^*(X;\bF_p)$ and $\boldk=H\bF_p$. It is not hard to see that  $\boldR$  is regular 
if and only if $H_*(\Omega X; \bF_p)$ is finite dimensional.

If $G$ is a finite $p$-group, $X=BG$ is $p$-complete and $\Omega BG \simeq G$, 
so that $\boldR$ is regular in this case. More generally, for $p$-complete spaces
$X$, regularity of $\boldR$ is equivalent to $X$ being the 
classifying space of a $p$-compact group in the sense \cite{DW} of Dwyer-Wilkerson.

\subsection{The Gorenstein condition.}
A commutative Noetherian local ring $R$ is Gorenstein if and only if $\Ext_{R}^*(k,R)$
is  one dimensional as a $k$-vector space. In the derived category we can restate
this as saying that the homology of the (right derived) Hom complex $\Hom_{R}(k,R)$ 
is equivalent to a suspension of $k$. This suggests the definition for ring spectra.

\begin{defn}
We say that $\boldR\lra \boldk$ is {\em Gorenstein} if there is an equivalence of $\boldR$-modules 
$\Hom_\boldR(\boldk,\boldR) \simeq \Sigma^a\boldk$ for some integer $a$.
\end{defn}

Evidently, $R$ is Gorenstein if and only if $\boldR=HR \lra Hk=\boldk$ is Gorenstein in 
the new sense.
However there is an interesting new phenomenon for spectra. We note that $\Hom_\boldR(\boldk,\boldR)$
admits a {\em right} action by the (derived) endomorphism ring $\cE = \HombR (\boldk,\boldk)$, 
whereas $\boldk$ is naturally a {\em left}
$\cE$-module. Thus if $\boldR$ is Gorenstein, $\boldk$ acquires new structure: that of a right
$\cE$-module. We want to say that $\boldR$ is {\em orientable} if this right action is 
trivial, but we must pause to give meaning to the notion of triviality.

\subsection{Orientability.}
We say that an $\boldR$-module $I$ is a {\em Matlis lift} of $\boldk$ if it is built from $\boldk$
using triangles and (arbitrary) coproducts and in addition $\HombR (\boldk,I) \simeq \boldk$ as
$\boldR$-modules. For example, if $R$ is a local ring, the injective hull $I(k)$ of $k$
is a Matlis lift of $k$.  If $\boldR$ is a $\boldk$-algebra then $I=\Hombk(\boldR,\boldk)$ is a Matlis
lift of $\boldk$ if it is built from $\boldk$.

In general there may be several Matlis lifts, or no Matlis lifts at all, 
but in many circumstances there is a preferred one. The above examples will 
cover the cases we consider, and we assume that a Matlis lift has been chosen. 
We use this to define what we mean by the trivial action of $\cE$ on $\boldk$
(i.e., namely the right action of $\cE$ on $\HombR (\boldk, I)$).

\begin{defn}
A Gorenstein ring spectrum $\boldR$ is {\em orientable} if $\HombR(\boldk,\boldR) \simeq \Sigma^a \HombR (\boldk,I)$
as right $\cE$-modules. 
\end{defn}

It turns out that for local rings $R$, there is a unique right $\cE$-module 
structure on $\boldk$, and hence every Gorenstein commutative ring is 
orientable as a ring spectrum, and the 
notion of orientability is irrelevant to classical commutative algebra.

However things are more interesting for the cochain algebra $\boldR=C^*(X;k)$
on a space $X$. By Poincar\'e duality, such a ring spectrum $\boldR$ is 
orientably Gorenstein if $X$ is a compact connected manifold orientable
over $k$. More generally, if $k=\Z /2^n$ the ring spectrum $\boldR$ is Gorenstein if $X$
is a compact connected manifold, and $\boldR$ is orientable if and only if 
the manifold $X$ is orientable over $k$. Indeed, right actions of $\cE $ on 
$\boldk$ correspond to group 
homomorphisms $\pi_1(X) \lra k^{\times}$, and the Gorenstein action of $\cE$
corresponds to the orientation homomorphism for the tangent bundle.

Similarly, the ring spectrum is orientably Gorenstein if $X=BG$ is 
the classifying space of a finite group, and more generally 
if $G$ is a compact Lie group with the property that the (adjoint) action of 
component group on the Lie algebra $T_eG$ is trivial over $k$.
More generally, if $k=\Z /2^n$, the ring spectrum is Gorenstein, 
and it is orientable if and only if the adjoint action is trivial 
over $k$. 

To exploit the Gorenstein condition to give structural information, we need to discuss
Morita equivalences. 

\subsection{Morita equivalences.}
Continuing to let $\cE=\HombR (\boldk,\boldk)$ denote the (derived) endomorphism
ring, we consider the relationship between the derived categories of left $\boldR$-modules 
and  of right $\cE$-modules. We have  the adjoint pair
$$\adjunction{T}{D(\modcE)}{D(\bRmod)}{E}$$
defined by 
$$T(X):=X \tensorcE \boldk \mbox{ and } E(M):=\HombR (\boldk,M).$$
It is easy to see that if $\boldk$ is finitely built from $\boldR$,  this gives an 
equivalence between the derived 
category of $\boldR$-modules built from $\boldk$
and the derived category of $\cE$-modules. However a variant will be more useful to us.
For this we say that $\boldk$ is {\em  proxy-small} if there is a finite $\boldR$-module
$K$ so that $\boldk$ is built from $K$ and $K$ is finitely built from $\boldk$. In 
the local algebra example with $\boldk=Hk$ (and $k$ a field), we may take 
$K$ to correspond to the Koszul complex and see that $\boldk$ is always proxy-small.
The variant  then reads as follows. 

\begin{lemma}  
Provided $\boldk$ is proxy-small, 
 the counit $TEM \lra M$ is an equivalence if $M$ is built from $\boldk$.
\qqed
\end{lemma}

\subsection{The local cohomology theorem.}
Now suppose that $\boldR \lra \boldk$ is proxy-small, and note that it is a 
formality that we can construct a $\boldk$-cellular 
approximation $\Gamma \boldR$ to $\boldR$. By definition, 
the $\boldk$-cellular approximation is an $\boldR$-module $\Gamma \boldR$ built from $\boldk$ with 
 a map $\Gamma \boldR \lra \boldR$ inducing an equivalence
$\HombR (\boldk,\boldR) \simeq \HombR (\boldk, \Gamma \boldR)$.
In the local algebra example,  we can take $\Gamma \boldR$ to be given by the 
stable Koszul complex. We may use a stable Koszul complex as a model for the 
$\boldk$-cellular approximation more generally, for example 
if $\pi_*(\boldR)$ is a Noetherian local ring 
with residue field $\pi_*(\boldk)$, and this shows its homotopy is calculated 
using local cohomology in that there is a spectral sequence 
$$
H^*_{\fm}(\pi_*(\boldR))\Rightarrow \pi_*(\Gamma \boldR).
$$

 We can then deduce a valuable duality property from the Gorenstein condition. 
Indeed, if $\boldR \lra \boldk$ is orientably Gorenstein, we have the equivalences
$$
E\Gamma \boldR =\HombR (\boldk,\Gamma \boldR) \simeq 
\HombR(\boldk,\boldR) \simeq \Sigma^a \HombR (\boldk, I) =E\Sigma^aI
$$ 
of right $\cE$-modules. Now applying Morita theory,   we conclude
$$\Gamma \boldR \simeq TE\Gamma \boldR \simeq TE \Sigma^aI \simeq \Sigma^a I.$$
For example if $\boldR$ is a $\boldk$-algebra with $\boldk =Hk$ for a field $k$, 
we can take $I =\Hombk (\boldR,\boldk)$ and conclude that
there is a spectral sequence
$$H^*_{\fm}(\pi_*(\boldR))\Rightarrow H^*(\Hombk (\boldR,\boldk))=\Homk (\pi_*\boldR,k).$$
In particular {\em if $\pi_*(\boldR)$ is Cohen-Macaulay}, this spectral sequence 
collapses to show {\em it is also Gorenstein}. In fact one may apply Grothendieck's
dual localization process to this spectral sequence and hence conclude that 
{\em in any case $\pi_*(\boldR)$ is generically Gorenstein} \cite{GL}.

For example we have seen that $C^*(BG)$ is regular if $G$ is a finite $p$-group; 
it follows that it is Gorenstein. Since $\pi_*(C^*(BG))=H^*(BG)$,  there is a spectral sequence
$$H^*_{\fm}(H^*(BG))\Rightarrow H_*(BG), $$
showing that the group cohomology ring $H^*(BG)$ has very special properties, 
such as being generically Gorenstein.

\bibliographystyle{amsalpha}

\end{document}